\documentstyle[12pt]{article}
\if@twoside  m
\else
\fi
\newcommand{\ie}{{\em i.e.}}
\newcommand{\eg}{{\em e.g.}}
\newcommand{\cf}{{\em cf. }}

\newcommand{\eps}{\varepsilon}
\newcommand{\DD}{{\cal D}}
\newcommand{\HH}{{\cal H}}
\newcommand{\OO}{{\cal O}}
\newcommand{\R}{I\!\!R}
\newcommand{\Z}{Z\!\!\!Z}
\newcommand{\Iff}{{\em iff $\:$}}
\newcommand{\re}{{\rm Re\,}}
\newcommand{\im}{{\rm Im\,}}
\newcommand{\QED}{\mbox{\rule[-1.5pt]{6pt}{10pt}}}
\newtheorem{claim}{Claim}[section]
\newtheorem{theorem}[claim]{Theorem}
\newtheorem{proposition}[claim]{Proposition}
\newtheorem{remark}[claim]{Remark}

\begin{document}
\thispagestyle{empty}
\title{\bf Contact interactions on graph superlattices}
\date{}
\author{Pavel Exner}
\maketitle
\begin{center}
Nuclear Physics Institute, Academy of Sciences \\25068 \v{R}e\v{z}
near Prague \\ and Doppler Institute, Czech Technical University, \\

B\v rehov\'a 7, 11519 Prague, \\ Czech Republic \\ {\em
exner@ujf.cas.cz}
\end{center}
\vspace{10mm}
\begin{quote}
{\small We consider a quantum mechanical particle living on a graph
and discuss the behaviour of its wavefunction at graph vertices. In
addition to the standard (or $\,\delta\,$ type) boundary conditions
with continuous wavefunctions, we investigate two types of a singular
coupling which are analogous to the $\,\delta'$ interaction and its
symmetrized version for particle on a line. We show that these
couplings can be used to model graph superlattices in which point
junctions are replaced by complicated geometric scatterers. We also
discuss the band spectra for rectangular lattices with the mentioned
couplings.  We show that they roughly correspond to their
Kronig--Penney analogues: the $\,\delta'$ lattices have bands whose
widths are asymptotically bounded and do not approach zero, while the
$\,\delta\,$ lattice {\em gap} widths are bounded. However, if the
lattice--spacing ratio is an irrational number badly approximable by
rationals, and the $\,\delta\,$ coupling constant is small enough,
the $\,\delta\,$ lattice has {\em no gaps} above the threshold of the
spectrum. On the other hand, infinitely many gaps emerge above a
critical value of the coupling constant; for almost all ratios this
value is zero.}
\end{quote}

\section{Introduction}
The problem of describing a quantum particle living on a graph is not
new in quantum mechanics; it appeared for the first time in early
fifties in connection with the free--electron model of organic
molecules \cite{RuS}. Writing down a Hamiltonian of such a particle
requires to check that the coupling between the wavefunctions at
branching point of the graph defines a self--adjoint operator, or in
physical terms, that the probability flow is preserved there. This is
conventionally achieved by demanding that the wavefunctions are
continuous at the junctions and satisfy the conditions
   \begin{equation} \label{delta bc}
\sum_j f'_j(x_m)\,=\, c_mf(x_m)\,,
   \end{equation}
where $\,m\,$ is the vertex number, the sum runs over all links
entering this vertex, $\,f(x_m)\,$ is the common value of the
functions $\,f_j\,$ there, and $\,c_m\,$ is a real parameter.

In recent years the interest to quantum mechanics on graphs has been
revived --- see, \eg, \cite{Ad,AL,AEL1,ARZ,BuT,ES1,GLR,GP} and
references therein --- in particular, as a reaction to the rapid
progress of fabrication techniques which allow us nowadays to produce
plenty of graph--like structures of a pure semiconductor material,
for which graph Hamiltonians represent a natural model. This posed
anew the question about physical plausibility of the boundary
conditions (\ref{delta bc}).

This problem has two basic aspects. The first of them concerns the
derivation of graph models from a more realistic description, in
which the configuration space consists of a system of coupled tubes.
This is still not the true system in which the tubes are complicated
many--body objects, but the crystallic structure of the semiconductor
material makes it a reasonable approximation.

Replacing a tube system by a graph of the same geometry, say, by the
tube axes, means a substantial simplification. The reason is that one
cannot {\em assume} that the wavefunction in a tube is independent of
the radial a azimuthal coordinates as the authors of \cite{GLR}
naively suggest. The tube Laplacian is specified by appropriate
boundary conditions, so there is a fully concrete system of
transverse eigenfunctions; most frequently the tube boundary is
assumed to be Dirichlet in which case none of these eigenfunctions is
constant. The graph approximation is generally expected to work in
situation where the tubes under consideration are thin enough so that
the transverse--mode eigenstates are well separated in energy and
their geometrically induced coupling coming from bending and
branching is weak.

The one--mode approximation for a single bent tube can be justified
\cite{DE,DES} but the problem is in no case a simple one. The case of
a branched tube is even harder because a typical branching region is
self--similar with respect to changes of the tube widths, so there is
no natural parameter to switch off the intermode coupling. In
general, we know neither the condition under which the graph
approximation works for a branched tube, nor the boundary conditions
which should model a particular branching geometry. However, we are
not going to discuss these problems in the present paper, apart from
some comments given in the concluding section.

The second aspect concerns intrinsic properties of the graph
Hamiltonians: one may ask what is the family of admissible operators
and which place is occupied in this class by those referring to the
conditions (\ref{delta bc}). This problem was solved in \cite{ES1}
where we showed how a general graph Hamiltonian can be constructed
using the von Neumann theory of self--adjoint extensions. However,
the operator family obtained in this way is large: even if we exclude
nonlocal interactions, \ie, we forbid particle hopping between
different branching points, each vertex of the graph is associated
with $\,n^2\,$ real parameters, where $\,n\,$ is the number of links
entering the junction.

Motivated by this we discussed in \cite{ES1} in detail several
subfamilies of such operators. The simplest situation occurs if the
domain of a graph Hamiltonian is required to consists of functions
which are continuous at each vertex. Then we arrive back at the
situation we started with, namely at the conditions (\ref{delta bc})
in which every junction is characterized by a single parameter.
Contrary to the claim made in \cite{GLR}, however, this parameter is
a real number and not a function of energy.

The above discussion shows that until derived within a
squeezing--tube approximation, the wavefunction continuity is just an
assumption which we may or may not adopt. If we decide to drop it,
the next more general class of graph Hamiltonians consists of those
which are {\em locally} permutation invariant at each vertex; we have
shown in \cite{ES1} that any junction is then described by a {\em
pair} of real parameters in such a way that
   \begin{equation} \label{permutation bc}
f_j(x_m)\,=\, A_mf'_j(x_m)\,+\,B_m\,\sum_{k\ne j} f'_k(x_m)\,,\quad
j\,=\,1,\dots,n_m\,,
   \end{equation}
where the indices have the same meaning as above; this condition
reduces to (\ref{delta bc}) for $\,A_m=B_m=:c_m^{-1}$. At the same
time, a junction can be described also by two one--parameter families
of boundary conditions, which represent singular limits of
(\ref{permutation bc}) when $\,A_m,B_m\to\pm\infty\,$ while
   \begin{equation} \label{singular limit}
C_m\,:=\,B_m-A_m\quad {\rm or} \quad D_m\,:=n_m\lbrack
A_m+(n_m\!-\!1)B_m\rbrack\,
   \end{equation}
is kept preserved; they are described by the boundary conditions
(\ref{C}) and (\ref{D}) below, respectively.

Although these couplings and their scattering properties were
discussed in detail in \cite{ES1}, their meaning remained somewhat
obscurred.  Our aim in this paper is to show that they represent in a
sense a counterpart to the coupling given by (\ref{delta bc}), and
that they generalize naturally the concept of $\,\delta'$ interaction
which has attracted attention recently in connection with spectral
properties of Wannier--Stark systems \cite{AEL1,AEL2,MS}. This will
be done in the next section; we are also going to show there that the
relation between $\,\delta'$ and geometric scatterers discovered in
\cite{AEL1} extends to vertices of any number of links.

After that we shall discuss the spectra of periodic rectangular
lattices in which each junction is described by one of the couplings
mentioned above. For the sake of simplicity, we restrict ourselves to
the planar case, however, the results have a straightforward
extension to higher dimensions. We shall demonstrate, in particular,
that the $\,\delta\,$ lattice spectrum depends substantially on the
ratio $\,\theta\,$ of the rectangle sides. It can even happen that
there are no gaps above the threshold of the spectrum; this is the
case if $\,\theta\,$ is badly approximable by rationals and the
coupling is ``weak'' enough. On the other hand, infinitely many gaps
exist for almost all $\,\theta\,$ but the gap pattern, as well as the
band pattern for a $\,\delta'_s\,$ lattice, is again irregular for an
irrational $\,\theta\,$. A summary of the results has been given in
\cite{KP}.

\setcounter{equation}{0}
\section{$\,\delta'$ and $\,\delta'_s\,$ interactions}

Throughout this section we will consider a single branching, hence we
may drop the index denoting the junction. For simplicity we will deal
with the graph $\,\Gamma_n\,$ consisting of $\,n\,$ halflines whose
endpoints are connected at a single point; as explained in
\cite{ES1}, the boundary conditions describing the coupling are local
and may be used for any $\,n$--link vertex.

The state Hilbert space for a spinless particle having $\,\Gamma_n\,$
as the configuration space is $\,\HH:= \bigoplus_{j=1}^n
L^2(0,\infty)\,$ and the free Hamiltonian acts at $\,f\equiv
\{f_j\}\,$ as $\,Hf=\{-f''_j\}\;$; its domain consists of all
functions whose components are $\,H^{2,2}$ on the halflines and
satisfy appropriate boundary conditions at the connection point which
relate the boundary values $\,f_j(0+)\,$ and $\,f'_j(0+)\,,
\;j=1,\dots,n\,$. For the sake of brevity we shall drop the arguments
as long as a single connection point is concerned.

\subsection{A warm--up: singular interactions on line}

Consider first the case of $\,\Gamma_2\,$ which is naturally
isomorphic a line with a single point interaction. In general, there
is a four--parameter family of such interactions ---
\cf\cite{ABD,Ca,CH,EG,Se1}. The best known among them are the
$\,\delta\,$ interactions for which the wavefunction is continuous at
the point supporting the interaction, while its derivative has a jump
proportional to the function value \cite[Sec.I.3]{AGHH}. Another
important class are the $\,\delta'$ interactions where the roles of
the functions and derivatives are reversed: the wavefunctions now
satisfy the conditions
   \begin{equation} \label{delta'}
f'_+=f'_-\,:=\,f'\,, \qquad f_+-f_-=\beta f'\,, \qquad \beta\in\R\,.
   \end{equation}
The name is somewhat misleading because in distinction to
$\,\delta\,$, the $\,\delta'$ interaction {\em cannot} be
approximated by Schr\"odinger operators with scaled potentials
\cite{AGHH,Se2}.  Instead, one can use approximations by families of
rank--one operators \cite{Se2} or velocity--dependent potentials
\cite{Ca}. Another way to understand the $\,\delta'$ interactions has
been suggested in \cite{AEL1}: replacing the line by a
``spiked--onion'' graph, \ie, cutting it into two halflines and
joining the loose ends by $\,N\,$ links of length $\,\ell\,$, one can
reproduce the high--energy scattering behaviour of the $\,\delta'$
interaction {\em up to a phase factor} in the limit when $\,N\to
\infty\,, \,\ell\to 0\,$ and $\,\beta:=N\ell\,$ is kept preserved.

In the two--parameter family of space--reversal invariant
interactions there is another subset with a similar property which
may be regarded as a {\em symmetrized version} of the $\,\delta'$
interaction; it is characterized by the boundary conditions
   \begin{equation} \label{delta's}
f'_++f'_-=0\,, \qquad f_++f_-=D f'_+\,, \qquad D\in\R\,.
   \end{equation}
It is easy to check that it has a bound state of energy $\,-4D^{-2}$
provided $\,D<0\,$ and the reflection and transmission amplitudes
for a plane wave of momentum $\,k\,$ are
   \begin{equation} \label{rt's}
r(k)\,=\, {-ikD\over 2-ikD}\,,\qquad t(k)\,=\, {-2\over 2-ikD}\,,
   \end{equation}
which up to signs coincides with both the $\,\delta'$ result (for
$\,\beta=D\,$) and the high--energy behaviour of the limiting
geometric scatterer.

\subsection{Vertices with $\,n\,$ links}

Let us look now how the above results extend to a quantum motion on
$\,\Gamma_n\,$ with $\,n\ge 3\,$. The coupling with continuous
wavefunctions is a natural analogue of the $\,\delta\,$ interaction,
so we shall use for it the same name. The boundary conditions
(\ref{delta bc}) can be rewritten as
   \begin{equation} \label{delta}
f_1=\cdots=f_n\,=:\,f\,, \qquad \sum_{j=1}^n\,f'_j\,=\,cf\,.
   \end{equation}
The first of the two classes mentioned in the introduction is for a
given $\,C\in\R\,$ characterized by the relations
   \begin{equation} \label{C}
\sum_{j=1}^n\,f'_j\,=\,0 \qquad {\rm and} \qquad
f_j\!-f_k+ C(f'_j\!-f'_k)\,=\,0\,, \quad j,k=1,\dots,n\;;
   \end{equation}
it is clear that just $\,n\,$ among these conditions are independent.
For $\,n=2\,$ this coincides with the requirements (\ref{delta'}) if
we put $\,\beta=-2C\;$; recall that both links are now positive real
halflines so one of the derivatives has to change sign.  In view of
this analogy {\em we shall refer to the coupling (\ref{C}) as to}
$\,\delta'$. It is sufficient to consider $\,C\ne 0\,$, because
otherwise the conditions (\ref{C}) reduce to (\ref{delta}) for
$\,c=0\,$.  On the other hand, the second exceptional class is given
by
   \begin{equation}
\label{D} f'_1=\cdots=f'_n\,=:\,f'\,, \qquad
\sum_{j=1}^n\,f_j\,=\,Df'
   \end{equation}
with $\,D\in\R\,$. In the two--link case this correspond to the
condition (\ref{delta's}), hence we shall speak about the
$\,\delta'_s\,$ {\em coupling.}

To justify the claim represented by the relation (\ref{singular
limit}), let us solve (\ref{permutation bc}) with respect to the
derivatives. This yields
   \begin{equation} \label{permutation bc2}
f'_j\,=\, {\DD_{n-2}\over(\!A-\!B)\DD_{n-1}}\,f_j\,-\,
{B\over(A\!-\!B)\DD_{n-1}}\, \sum_{k\ne j}\,f_k\,,
   \end{equation}
where $\,\DD_m:= A+mB\,$. In the first limit of (\ref{singular
limit}) we arrive then at
$$
Cf'_j\,=\, {1-n\over n}\,f_j\,+\, {1\over n}\, \sum_{k\ne j}\,f_k\,,
$$
which is equivalent to (\ref{C}), while the other limit gives
$$
f'_j\,=\, {1\over D}\, \sum_{k=1}^n\, f_k\,,  \qquad j=1,\dots,n\,,
$$
\ie, the condition (\ref{D}).

Each of the three operators on $\,\Gamma_n\,$ described above can
exhibit bound states with eigenfunctions localized around the
junction --- this happens \Iff the corresponding parameter,
$\,c,\;-C\,$ or $\,D\,$, respectively, is negative (we keep the
notation of \cite{ES1}; it would be more natural in the present
context to change the sign of $\,C\,$). Since any solution
to the free Schr\"odinger equation is of the form $\,\{\alpha_j\,
e^{-\kappa x_j}\}\,$, it is sufficient to substitute to the boundary
conditions (\ref{delta})--(\ref{D}) to check that the corresponding
eigenvalues are
$$
\epsilon(c)\,=\,-\,\left(c\over n\right)^2\,, \quad
\epsilon(C)\,=\,-\,C^{-2} \qquad {\rm and} \qquad
\epsilon(D)\,=\,-\,\left(n\over D\right)^2
$$
in the three cases; the $\,\delta\,$ and $\,\delta'_s\,$ eigenvalues
are simple, while the $\,\delta'$ bound state has multiplicity
$\,n\!-\!1\,$. Also the scattering properties of these junctions are
easily found --- see \cite{ES1}.

\subsection{Geometric--scatterer junctions}

Now we want to show that the ``spiked--onion'' argument mentioned in
Section~2.1 can be naturally extended to $\,\delta'$ and
$\,\delta'_s\,$ couplings on $\,\Gamma_n\,,\, n\ge 3\,$. Let us
replace the latter by the graph $\,\Gamma_n(N,\ell)\,$ sketched on
Figure~1: every pair of halfline endpoints is connected by $\,N\,$
links of length $\,\ell\;$; the corresponding variables will run
through the interval $\,[-\ell/2,\ell/2]\,$. We shall assume that the
coupling at each graph node is given by the condition (\ref{delta})
with the same parameter $\,c\,$.

%

Using the permutation symmetry of the Hamiltonian, we can write a
general scattering--solution Ansatz in the following form:
   \begin{description}
   \item{  (i)} $\;e^{-ikx}+r\,e^{ikx}\;$ at a chosen external
link (halfline),
   \item{ (ii)} $\;t\,e^{ikx}\;$ at the other $\,n\!-\!1\,$
halflines,
   \item{(iii)} $\;\alpha\,e^{ikx}+\beta\,e^{-ikx}\;$ at the
$\,2N\,$ connecting links coupled to the ``incident'' halfline,
   \item{ (iv)} $\;\delta\,\cos kx\;$ at the remaining
$\,N\left\lbrack {n\choose 2}-2 \right\rbrack\,$ connecting links.
   \end{description}
The condition (\ref{delta}) now yields the following system of
equations
   \begin{eqnarray} \label{matching}
1+r &\!=\!& \alpha\bar\eta+ \beta\eta  \nonumber \\ \nonumber \\
t\,=\, \alpha\eta+ \beta\bar\eta &\!=\!& \delta \cos{k\ell\over 2}
\nonumber \\ \\
r-1+N(n\!-\!1)(\alpha\bar\eta- \beta\eta) &\!=\!& \gamma(1+r)
\nonumber \\ \nonumber \\
t-N(\alpha\eta- \beta\bar\eta) -iN(n\!-\!2)\delta\, \sin{k\ell\over
2} &\!=\!& \gamma t\,,  \phantom{AAAAAAAAAA} \nonumber
   \end{eqnarray}
where we have denoted
$$
\eta\,:=\, e^{ik\ell/2}\,, \qquad \gamma\,:=\,{c\over ik}\,.
$$
The first three equations coming from the continuity requirement are
solved by
$$
\alpha\,=\, {i\over 2}\, {\bar\eta(1\!+\!r)-\eta t\over \sin
k\ell}\,, \qquad \beta\,=\, {i\over 2}\, {\bar\eta
t-\eta(1\!+\!r)\over \sin k\ell}\,, \qquad \delta\,=\, {t\over
\cos{k\ell\over 2}}\,.
$$
Next we substitute these values into the remaining two of the
equations (\ref{matching}) and use the identity
$$
\cot 2x- (n\!-\!2)\tan x\,=\, (n\!-\!1)\cot 2x- (n\!-\!2)\csc 2x\;;
$$
introducing
$$
P\,\equiv \,P(k)\,:=\, 1-\gamma+iN(n\!-\!1)\cot k\ell\,,
\qquad Q\,\equiv Q(k)\,:=\, {iN\over \sin k\ell}\,,
$$
we can rewrite the resulting system as
$$
rP-(n\!-\!1)Qt\,=\, \bar P\,, \qquad -rQ+ t[P-(n\!-\!2)Q]\,=\, Q\,.
$$
Since $\,\re P=1\,$, the sought reflection and transmission
amplitudes are given by
   \begin{eqnarray} \label{rt}
r(k) &\!=\!& {|P|^2\!-(n\!-\!2)\bar PQ+ (n\!-\!1)Q^2\over
P^2\!-(n\!-\!2) PQ+ (n\!-\!1)Q^2}\,, \nonumber \\ \\
t(k) &\!=\!& {2Q\over
P^2\!-(n\!-\!2) PQ+ (n\!-\!1)Q^2}\,. \nonumber
   \end{eqnarray}
It is straightforward to check the identity
$$
|P^2\!-(n\!-\!2) PQ+ (n\!-\!1)Q^2|^2- ||P|^2\!-(n\!-\!2)\bar PQ+
(n\!-\!1)Q^2|^2 \,=\,4(n\!-\!1)|Q|^2\,,
$$
so the S--matrix is unitary,
   \begin{equation} \label{unitarity}
|r(k)|^2+(n\!-\!1)|t(k)|^2 \,=\,1\,.
   \end{equation}
Consider now the limit of increasingly complicated scatterers; in
other words, we put $\,\ell:={\tau\over N}\,$ and let
$\,N\to\infty\,$. The expressions (\ref{rt}) then can be for large
$\,N\,$ rewritten as
   \begin{eqnarray} \label{limiting rt}
r(k) &\!=\!& {n-2-\,{nc\over ik}\,-\,{n\choose 2} ik\tau\over
-n+\,{nc\over ik}\,+\,{n\choose 2} ik\tau}\,+\,\OO(N^{-1})\,,
\nonumber \\ \\
t(k) &\!=\!& {-2\over-n+\,{nc\over ik}\,+\,{n\choose 2}
ik\tau}\,+\,\OO(N^{-1})\,, \nonumber
   \end{eqnarray}
from where the limits immediately follow. At large energies, in
particular, the terms containing $\,c\,$ can be neglected and the
amplitudes (\ref{limiting rt}) behave as
   \begin{equation} \label{high-energy rt}
r(k) \,\approx\, {n-2-\,{n\choose 2} ik\tau\over
-n+\,{n\choose 2} ik\tau}\,, \qquad
t(k) \,\approx\, {-2\over-n+\,{n\choose 2} ik\tau}\;;
   \end{equation}
if $\,c=0\,$, of course, the right sides give the expressions of the
limiting amplitudes for any $\,k\,$.

Let us now compare this result to the S--matrices corresponding to
our two singular couplings which have been computed in \cite{ES1}. In
the $\,\delta'$ case the reflection and transmission amplitudes are
   \begin{equation} \label{rt C}
r(k) \,=\, {2-n+inkC\over n+inkC}\,, \qquad
t(k) \,=\, {2\over n+inkC}\;;
   \end{equation}
choosing $\,\tau:=\, -\,{2C\over n-1}\,$, we get at high energies the
same $\,t(k)\,$, while the reflection amplitude differs by the phase
factor $\,-2\arg(n-2+inkC)\,$ which goes to $\,\pi\,$ as
$\,k\to\infty\,$. As for the $\,\delta'_s\,$ coupling, we have
   \begin{equation} \label{rt D}
r(k) \,=\, {n-2-ikD\over n-ikD}\,, \qquad
t(k) \,=\, {-2\over n-ikD}\,,
   \end{equation}
so one has to put $\,\tau:=\,{2D\over n(n-1)}\,$ to recover the
S--matrix elements (\ref{high-energy rt}) up to a sign, which is
switched for all of them.

Hence both the singular couplings reproduce the high--energy
behaviour of the limiting geometric scatterer up to a phase factor.
They are the only ones with this property in the class
(\ref{permutation bc}). To see this, consider the reflection
amplitude corresponding to (\ref{permutation bc}). Leaving aside the
case $\,A=B\,$ corresponding to (\ref{delta}) when we have
$\,\lim_{k\to\infty} r(k)= {2-n\over n}\,$, its high--energy
behaviour is
$$
r(k)\,=\, {n-2+ik\left(1-{A\over B}\right)\DD_{n-1}\over
-n+2\left(1-{A\over B}\right) +ik\left(1-{A\over B}\right)
\DD_{n-1}}\,+\,\OO(k^{-2})\,,
$$
which certainly differs from (\ref{high-energy rt}) by more than a
phase factor. On the other hand, all these couplings (with
exception of the case $\,A=B\,$) represent the effective {\em
Neumann} decoupling at high energies, $\,\lim_{k\to\infty} r(k)=
1\,$, while the geometric scatterer mimicks rather the {\em
Dirichlet} decoupling, $\,\lim_{k\to\infty} r(k)= -1\,$.

\setcounter{equation}{0}
\section{Lattices with a singular coupling}

Consider now the case which has attracted some attention recently
\cite{AL,GLR} as a model of quantum--wire superlattices. We shall
assume that the graph in question is a rectangular lattice with the
spacings $\,\ell_1,\,\ell_2\,$ in the $\,x\,$ and $\,y\,$ direction,
respectively (\cf Figure~2). In addition, we suppose that each graph
vertex is endowed with {\em the same} coupling of one of the above
described types. We restrict ourselves to the planar situation just
for the sake of simplicity; the band conditions obtained below and
the method of their solution have a straightforward extension to
higher dimensions.

%

\subsection{The Bloch analysis}

To find the band spectrum of such a lattice, we start from a natural
Ansatz for the Bloch solutions: we choose
   \begin{eqnarray} \label{Bloch Ansatz}
f_m(x) &\!=\!& e^{im\theta_2\ell_2} \left( a_n e^{ikx}+b_n e^{-ikx}
\right) \qquad \dots \qquad x\in(n\ell_1,(n\!+\!1)\ell_1) \nonumber
\\ \\ g_n(y) &\!=\!& e^{in\theta_1\ell_1} \left( c_m e^{iky}+d_m
e^{-iky} \right) \qquad \dots \qquad y\in(m\ell_2,(m\!+\!1)\ell_2)
\nonumber
   \end{eqnarray}
The coefficients on neighbouring links are related by
   \begin{eqnarray} \label{Bloch condition}
a_{n+1}\xi\,e^{ikx}+ b_{n+1}\bar\xi\,e^{-ikx} &\!=\!&
\sigma(a_{n}e^{ikx}+ b_{n}e^{-ikx}) \,, \nonumber \\ \\
c_{m+1}\eta\,e^{iky}+ d_{m+1}\bar\eta\,e^{-iky} &\!=\!&
\tau(b_{m}e^{iky}+ d_{m}e^{-iky}) \,, \nonumber
   \end{eqnarray}
where we have denoted for the sake of brevity
$$
\sigma\,:=\,e^{i\theta_1\ell_1}\,,\quad
\tau\,:=\,e^{i\theta_2\ell_2}\,,\quad \xi\,:=\,e^{ik\ell_1}\,,\quad
\eta\,:=\,e^{ik\ell_2}\;;
$$
this conditions allow us to compute easily the needed boundary
values.

Let us begin with the $\,\delta\,$ coupling, where the relations
(\ref{delta}) now read
   \begin{eqnarray} \label{delta lattice}
f_m(n\ell_1\!+\!0)\,=\, f_m(n\ell_1\!-\!0)\,=\,
g_n(m\ell_2\!+\!0)\,=\, g_n(n\ell_2\!-\!0) &\!=:\!& F_{mn}  \,,
\nonumber \\ \\ f'_m(n\ell_1\!+\!0)- f'_m(n\ell_1\!-\!0)+
g'_n(m\ell_2\!+\!0)- g'_n(n\ell_2\!-\!0)
&\!=\!& cF_{mn}\,. \nonumber
   \end{eqnarray}
Substituting the boundary values, we get a homogeneous system a four
independent equations for the coefficients, which has a solution
provided
$$
-2ik\bar\tau(\eta-\bar\eta)[1-\bar\sigma(\xi+\bar\xi)+ \bar\sigma^2]
\,+\, c\bar\sigma\bar\tau(\xi-\bar\xi)(\eta-\bar\eta)\,-\,
2ik\bar\sigma (\xi-\bar\xi)[1-\bar\tau(\eta+\bar\eta)+ \bar\tau^2]
\,=\,0\,.
$$
Returning to the original quantities, we can cast it into the form
   \begin{equation} \label{delta band condition}
{\cos\theta_1\ell_1-\cos k\ell_1\over \sin k\ell_1}\,+\,
{\cos\theta_2\ell_2-\cos k\ell_2\over \sin k\ell_2}\,-\,{c\over 2k}
\,=\,0 \;;
   \end{equation}
this is the result obtained in \cite{GLR}, and in \cite{AL} for the
particular case $\,c=0\,$.

If we assume instead that the lattice links are coupled by the
$\,\delta'_s\,$ interaction, the requirement (\ref{delta lattice}) is
replaced by
   \begin{eqnarray} \label{delta's lattice}
f'_m(n\ell_1\!+\!0)\,=\,-f'_m(n\ell_1\!-\!0)\,=\,
g'_n(m\ell_2\!+\!0)\,=\, -g'_n(n\ell_2\!-\!0) &\!=:\!& G_{mn} \,,
\phantom{AAAA} \nonumber \\ \\ f_m(n\ell_1\!+\!0)+
f_m(n\ell_1\!-\!0)+ g_n(m\ell_2\!+\!0)+ g_n(n\ell_2\!-\!0)
&\!=\!& DG_{mn}\;; \nonumber
   \end{eqnarray}
solving it in the same way we arrive at the condition
   \begin{equation} \label{delta's band condition}
{\cos\theta_1\ell_1+\cos k\ell_1\over \sin k\ell_1}\,+\,
{\cos\theta_2\ell_2+\cos k\ell_2\over \sin k\ell_2}\,-\,{Dk\over 2}
\,=\,0\,.
   \end{equation}
Finally, in the $\,\delta'$ case the coupling conditions read
   \begin{eqnarray} \label{delta' lattice}
f'_m(n\ell_1\!+\!0) -f'_m(n\ell_1\!-\!0) +g'_n(m\ell_2\!+\!0)
-g'_n(m\ell_2\!-\!0) &\!=\!& 0\,,  \nonumber \\ \nonumber \\
f_m(n\ell_1\!+\!0)- g_n(m\ell_2\!+\!0)+ C(f'_m(n\ell_1\!+\!0)-
g'_n(m\ell_2\!+\!0))   &\!=\!& 0\,, \nonumber \\ \\
g_n(m\ell_2\!+\!0)- f_m(n\ell_1\!-\!0) +C(g'_n(m\ell_2\!+\!0)+
f'_m(n\ell_1\!-\!0)) &\!=\!& 0\,, \nonumber \\ \nonumber \\
f_m(n\ell_1\!-\!0)- g_n(m\ell_2\!-\!0)+ C(-f'_m(n\ell_1\!-\!0)+
g'_n(m\ell_2\!-\!0))
&\!=\!& 0\,. \nonumber
   \end{eqnarray}
The solvability condition is now slightly more complicated; it can be
written as
   \begin{equation} \label{delta' band condition}
\sum_{j=1}^2\, {\cos\theta_j\ell_j- \re((1\!-\!ikC)e^{ik\ell_j}) \over
\im((1\!+\!ikC)^{-2}e^{ik\ell_j}) } \,=\,0\,.
   \end{equation}
At high energies, however, this relation simplifies: up to
$\,\OO(k^{-1})\,$ terms it acquires the form
   \begin{equation} \label{he delta' band condition}
{\cos\theta_1\ell_1-\cos k\ell_1\over \sin k\ell_1}\,+\,
{\cos\theta_2\ell_2-\cos k\ell_2\over \sin k\ell_2}\,-\,2Ck \,=\,0\,.
   \end{equation}
For a single junction considered in Section~2.3, a $\,\delta'_s\,$
coupling with the parameter $\,D\,$ has, up to a phase factor, the
same asymptotic behaviour as $\,\delta'$ for $\,C=-D/n\,$. Using
this substitution for $\,n=4\,$ we arrive at the relation
(\ref{delta's band condition}) with the parameters $\,\theta_j\ell_j\,$
replaced by $\,\theta_j\ell_j+\pi\,$. Hence the band spectra of the
$\,\delta'$ and $\,\delta'_s\,$ lattices behave alike at high
energies; we shall discuss below only the second case which is
simpler.

\subsection{$\,\delta\,$ lattice spectra}

The condition (\ref{delta band condition}) yields no restriction on
$\,k\,$ if $\,c=0\,$ as the authors of \cite{AL} pointed out; the
spectrum covers the interval $\,[0,\infty)\,$ and in the isotropic
case, $\,\ell_1=\ell_2\,$, one can write the energy $\,\epsilon(k^2)
:=k^2\,$ in terms of the Bloch parameters (quasimomentum components)
$\,\theta_1,\,\theta_2\,$ explicitly.

This is no longer true if the coupling constant is nonzero.
Nevertheless, one can say a lot about the spectrum determined by
(\ref{delta band condition}). To solve this condition, let us rewrite
it as
   \begin{equation} \label{delta band condition 2}
{c\over 2k}\,=\, F(k;v_1,v_2)\,:=\, \sum_{j=1}^2\, {v_j-\cos
k\ell_j\over \sin k\ell_j}\,,
   \end{equation}
where $\,v_j:=\cos\theta_j\ell_j\,$. Since these parameters run
through the interval $\,[-1,1]\,$, we find easily
   \begin{eqnarray*}
-\cot{k\ell_j\over 2} \,\le\, {v_j-\cos k\ell_j\over \sin k\ell_j}
\,\le\, \tan{k\ell_j\over 2} & \quad\dots\quad & \sin k\ell_j>0\,, \\
\\ \tan{k\ell_j\over 2} \,\le\, {v_j-\cos k\ell_j\over \sin k\ell_j}
\,\le\, -\cot{k\ell_j\over 2} & \quad\dots\quad & \sin k\ell_j>0\,.
   \end{eqnarray*}
{}From here we obtain the extremum values of $\,F(k;\cdot,\cdot)\,$ for a
fixed $\,k\,$:
   \begin{eqnarray} \label{delta extrema}
F_+(k) &\!:=\!& \max_{v_j\in[-1,1]} F(k;v_1,v_2)\,=\, \sum_{j=1}^2\,
\tan\left( {k\ell_j\over 2}\,-\,{\pi\over 2} \left\lbrack
{k\ell_j\over\pi} \right\rbrack \right) \,,  \nonumber \\ \\
F_-(k) &\!:=\!& \min_{v_j\in[-1,1]} F(k;v_1,v_2)\,=\,
-\,\sum_{j=1}^2\, \cot\left( {k\ell_j\over 2}\,-\,{\pi\over 2}
\left\lbrack {k\ell_j\over\pi} \right\rbrack \right) \,, \nonumber
   \end{eqnarray}
where the symbol $\,[\cdot]\,$ denotes the integer part; at the
singular points the value is understood as the limit from the right
or left, respectively. By definition we have  $\,\pm F_{\pm}(k)\ge
0\,$, hence the {\em gap} condition can be expressed as
   \begin{equation} \label{delta gap}
\pm\,{c\over 2k}\,>\, \pm F_{\pm}(k)
   \end{equation}
for $\,\pm c>0\,$, respectively. This applies to positive enegies,
on the negative halfline one compares instead $\,{c\over 2\kappa}\,$,
where $\,\kappa:=-ik\,$, with the extremum values
   \begin{equation} \label{delta extrema-}
F_+(\kappa) \,:=\, -\,\sum_{j=1}^2\,\tanh\left( {\kappa\ell_j\over 2}
\right)\,, \qquad F_-(\kappa) \,:=\, -\,\sum_{j=1}^2\,\coth\left(
{\kappa\ell_j\over 2} \right)\;;
   \end{equation}
since both of them are negative, there is obviously no negative
spectrum for $\,c\ge 0\,$. Notice that the band condition for a
one--dimensional array of $\,\delta\,$ interactions can be rewritten
in the same form with a single term on the right sides of the
relations (\ref{delta extrema}) and (\ref{delta extrema-}); it is
easy to reproduce from here the spectrum of the Kronig--Penney model
\cite[Sec.III.2]{AGHH}.

Let us collect first some simple properties of the spectrum which
follow directly from the condition (\ref{delta gap}) and
its negative--energy counterpart. We introduce
   \begin{equation} \label{parameters}
\ell\,:=\,\sqrt{\ell_1\ell_2}\,,\qquad \theta\,:=\,
{\ell_2\over\ell_1} \qquad {\rm and} \qquad
L\,:=\,\max(\ell_1,\ell_2)\,,
   \end{equation}
so
$$
\ell_1\,=\,\ell\theta^{-1/2}\,,\qquad
\ell_2\,=\,\ell\theta^{1/2}\,,\qquad
\min(\ell_1^{-1},\ell_2^{-1})\,=\, L^{-1}\,,
$$
and denote by $\,\sigma(\ell,\theta,c)\,$ the spectrum of the
corresponding $\,\delta\,$ lattice Hamiltonian.

   \begin{proposition} \label{delta basic}
$\;$ (a) The spectrum has a band structure, $\,\sigma(\ell,\theta,c)=
\bigcup_{r=1}^N [\alpha_r,\beta_r]\,$ for some $\,N\ge 1\,$, where
$\,\alpha_r<\beta_r<\alpha_{r+1}\,$.
\vspace{.8em}

\noindent
(b) $\;\sigma(\ell,\theta,0)= [0,\infty)\,$.
\vspace{.8em}

\noindent
(c) If $\,c>0\,$, each $\,\beta_r\,$ is of the form $\,\left({\pi
n\over\lambda} \right)^2\,$ for $\,\lambda=\ell_1\,$ or $\,\ell_2\,$,
and some $\,n\in\Z\,$. Similarly, $\,\alpha_r=\left({\pi
m\over\lambda} \right)^2\,$ with $\,m\in\Z\,$ for $\,c<0\,$ and
$\,r\ge 2\,$.
\vspace{.8em}

\noindent
(d) $\;\pm\alpha_1>0\;$ holds iff $\;\pm c>0\,$.
\vspace{.8em}

\noindent
(e) For $\,c< -4\ell^{-1}(\theta^{1/2}\!+\theta^{-1/2})\,$,
$\;\beta_1<0\,$ and $\;\alpha_2= \left(\pi\over L\right)^2\,$.
\vspace{.8em}

\noindent
(f) $\,\sigma(\ell,\theta,c')\cap\R_+ \,\subset\,
\sigma(\ell,\theta,c)\cap\R_+\,$ holds for $\,|c'|>|c|\,$.
\vspace{.8em}

\noindent
(g) Each gap is contained in the intersection of a pair of gaps of
the Kronig--Penney model with the coupling constant $\,c\,$ and
spacings $\,\ell_1\,$ and $\,\ell_2\,$.
\vspace{.8em}

\noindent
(h) All gaps above the threshold are finite. If there is an infinite
number of them, their widths are asymptotically bounded,
   \begin{equation} \label{gap width bound}
\alpha_{r+1}-\beta_r\,<\, 2|c|(\ell_1\!+\! \ell_2)^{-1}
+\,\OO(r^{-1})\,.
   \end{equation}
   \end{proposition}
\vspace{1em}

The property (g) allows us to localize spectral gaps by those of the
Kronig--Penney model. By negation, it illustrates that transport
properties of the lattice are better than a combination of its
$\,x\,$ and $\,y\,$ projections. If a given energy is contained in a
band in one of the directions, then by (\ref{delta band condition})
it is trivially also in a band of the lattice Hamiltonian. The
converse is not true, of course: the condition may be satisfied even
if none of the factors can be annulated separately. Of course,
different solutions correspond to different quasimomenta and
different directions in which the particle is able to ``dribble"
through the lattice.

However, the most interesting property of the spectrum is its
irregular dependence on $\,\theta\,$ coming from the existence of
competing periods in (\ref{delta extrema}). To formulate the results,
we have to recall some notions from the number theory \cite{HW,Sch}.
Irrational numbers can be classified with respect to how ``fast''
they can be approximated by rationals. In particular, such a number
is called {\em badly approximable} if there is a $\,\delta>0\,$ such
that
   \begin{equation} \label{bad}
\left|\, \theta\,-\,{p\over q}\right|\,>\,{\delta\over q^2}\,.
   \end{equation}
This is a non--empty subset in the family of all Diophantine numbers;
for instance, it contains all algebraic numbers of degree $\,2\,$,
\ie, irrational solutions of a quadratic equation with rational
coefficients. On the other hand, the Lebesgue measure of this set
is zero.

One can also write $\,\theta\,$ as a continued fraction
$\,[a_0,a_1,\dots]\,$ with integer coefficients; such a
representation is unique and provides a natural way to gauge the
approximation properties. The faster the $\,a_n$'s grow, the better
is $\,\theta\,$ approximated by the truncated fractions; the worst
irrational from this point of view is the golden mean $\,{1\over
2}\,(1+\sqrt 5) =[1,1,1,\dots]\,$. A counterpart to badly
approximable numbers is the class of irrationals with an
unbounded sequence of coefficients,
$$
\lim\!\sup_{\hspace{-1.5em} n\to\infty} a_n\,=\,\infty\,,
$$
which has the following property: there are sequences
$\,\{m_r\}_{r=1}^{\infty}\,, \; \{n_r\}_{r=1}^{\infty}\,$ of pairwise
relatively prime integers such that
   \begin{equation} \label{Last}
\lim_{r\to\infty} n_r^2\left|\, \theta\,-\,{m_r\over
n_r}\right|\,=\,0\,.
   \end{equation}
These numbers, which may be called {\em Last admissible} \cite{La},
have full Lebesgue measure. It is sometimes convenient to express
these approximations using number--theory symbols: the fractional
part $\,\{x\}:=x-[x]\,$, and $\,\|x\|:=\min(\{x\},1\!-\!\{x\})\,$.
Since the number in question may be indexed, to distiguish the
fractional part from a sequence, we shall always specify the range of
indices in the latter case.

   \begin{theorem} \label{delta number}
$\;$ (a) If $\,\theta\,$ is rational, $\,\sigma(\ell,\theta,c)\,$ has
infinitely many gaps for any nonzero coupling constant $\,c\,$.
\vspace{.8em}

\noindent
(b) For badly a approximable $\,\theta\,$ there is $\,c_0>0\,$ such
that for $\,|c|<c_0\,$ the spectrum has no gaps above the threshold,
$\,\beta_1=\infty\,$.
\vspace{.8em}

\noindent
(c) $\;\sigma(\ell,\theta,c)\,$ has infinitely many gaps for any
$\,\theta\,$ provided $\,|c|L> 5^{-1/2}\pi^2\,$.
\vspace{.8em}

\noindent
(d) If $\,\theta\,$ is Last admissible, there are infinitely many
gaps for any $\,c\ne 0\,$.
\vspace{.8em}
   \end{theorem}

\noindent
{\em Proof:} We shall suppose that $\,c>0\;$; the argument for a
negative coupling constant is analogous. If $\,\theta= {p\over q}\,$,
the function $\,F_+\,$ has infinitely many zeros in $\,\R_+\,$
without accumulation points, so (a) follows. On the other hand,
$\,F_+(k)>0\,$ for any $\,k>0\,$ if $\,\theta\,$ is irrational, hence
we have to investigate local minima of this function. They occur at
$\,k_n:=\,{\pi n\over\ell_1}\,$ and $\,\tilde k_m:=\,{\pi
m\over\ell_2}\,$ with $\,n,m\in\Z_+\,$, and the corresponding values
are
$$
F_+(k_n)\,=\, \tan\left({\pi\over 2}\{n\theta\}\right)\,, \qquad
F_+(\tilde k_m)\,=\, \tan\left({\pi\over 2}\{m\theta^{-1}\}\right)\,.
$$
If $\,\theta\,$ is badly approximable, the condition (\ref{bad})
yields
$$
F_+(k_n)\,>\, {\pi\over 2}\{n\theta\}\,\ge\, {\pi\over 2}\|n\theta\|
\,>\, {\pi\delta\over 2n}\,.
$$
Since $\,\theta^{-1}$ is also badly approximable, we get
$\,F_+(\tilde k_m)>\, {\pi\delta\over 2m}\;$; if different constants
correspond to $\,\theta,\,\theta^{-1}$, we call $\,\delta\,$ the
smaller of the two. This implies that for $\,c\,$ small enough
$$
F_+(k_n)\,>\, {c\over 2k_n}\,,\qquad F_+(\tilde k_m)\,>\, {c\over
2\tilde k_m}\,,
$$
holds for all $\,n,m\in\Z_+\,$, \ie, the assertion (b).

By the Hurwitz extension of Dirichlet's theorem \cite[Sec.II.1]{Sch}
one can find to any irrational $\,\theta\,$ sequences
$\,\{n_r\}_{r=1}^{\infty},\; \{m_r\}_{r=1}^{\infty}\,$ in such a way
that $\,|n_r\theta-m_r|< 5^{-1/2} n_r^{-1}$. Moreover, these
approximations can be constructed explicitly in terms of truncated
continued fractions \cite[Chap.10]{HW}. Choosing the truncations of even
lengths (without
the integer part), we obtain a sequence $\,\{{m_r\over
n_r}\}_{r=1}^{\infty}\,$ approaching $\,\theta\,$ from below. In that
case $\,\{n_r\theta\} \to 0\,$, and we get the estimate
$$
F_+(k_{n_r})\,<\,  {\pi\over 2}\,(1\!+\!\eps)\{n_r\theta\} \,<\,
{\pi(1\!+\!\eps) \over 2\sqrt{5}\, n_r}
$$
for any $\,\eps>0\,$ and $\,r\,$ large enough. In the same way one
can approximate $\,\theta^{-1}$. By (\ref{delta gap}) we find that
there are infinitely many gaps if $\,c> 5^{-1/2}\pi^2 L^{-1}
(1\!+\!\eps)\,$ for any $\,\eps>0\,$, \ie, the assertion (c).

Let finally $\,\theta\,$ be a Last admissible number. Without loss of
generality we may suppose that
$$
\lim\!\sup_{\hspace{-1.5em} n\to\infty} a_{2n}\,=\,\infty
$$
holds for its continued--fraction representation; otherwise we use
instead
$$
\theta^{-1}\,=\, \left\lbrace\; \begin{array}{lll} [0,a_0,a_1,\dots]
& \qquad\dots\qquad & a_0\ne 0 \\ \\ \lbrack a_1,a_2,a_3,\dots\rbrack
& \qquad\dots\qquad & a_0= 0 \end{array} \right.
$$
By (\ref{Last}) there is a sequence $\,\{n_r\}_{r=1}^{\infty}\,$ such
that $\,n_r\{n_r\theta\}\to 0+\,$. This means that for all large
enough $\,r\,$ we have
$$
n_r F_+(k_{n_r})\,=\, n_r\, \tan\left({\pi\over
2}\{n_r\theta\}\right) \,<\, \pi n_r\,\{n_r\theta\} \,\to\,0\,.
$$
Hence $\,k_{n_r} F_+(k_{n_r})\to 0\,$ too, so there are infinitely
many values of $\,k\,$ accumulating at infinity for which
$\,kF_+(k)<\, {c\over 2}\,$. \quad \QED

   \begin{remark}
{\rm The critical value of $\,\pi^2/\sqrt{5}= 4.414\dots\,$ in the
part (c) cannot be pushed down  --- \cf\cite[Sec.II.1]{Sch}. The
bound is saturated for the golden mean, $\,\theta={1\over
2}\,(1+\sqrt 5)\,$, and moreover, there are {\em no} gaps if
$\,|c|L\,$ is below the critical value in this case.}
   \end{remark}

\subsection{$\,\delta'_s\,$ lattice spectra}

Let us turn to lattices with a $\,\delta'_s\,$ coupling. The
relation (\ref{delta's band condition}) can be rewritten as
   \begin{equation} \label{delta's band condition 2}
{Dk\over 2}\,=\, -F(k;-v_1,-v_2)\,,
   \end{equation}
so the {\em bands} are now determined by the inequalities
   \begin{equation} \label{delta'_s band}
\mp F_{\mp}(k) \,\ge\,\pm\,{Dk\over 2}\,,\qquad \pm D>0\,.
   \end{equation}
This concerns positive energies; on the negative halfline we have
instead the condition
   \begin{equation} \label{delta'_s band-}
F_+(\kappa)\,\ge\, {D\kappa\over 2}\,\ge\, F_-(\kappa)
   \end{equation}
with $\,F_{\pm}(\kappa)\,$ given by (\ref{delta extrema-}); the
change in sign is due to the fact that $\,k\,$ is now in numerator on
the left side of (\ref{delta's band condition 2}). Replacing
$\,F_{\pm}(k)\,$ by similar expressions containing a single term, we
obtain in this way the band spectrum of the one--dimensional array of
$\,\delta'_s\,$ interactions.

   \begin{remark} \label{delta' rem}
{\rm Notice that the latter coincides with that of the $\,\delta'$
array of the same parameters \cite[Sec.III.3]{AGHH}, because the
corresponding transfer matrices differ just by sign; in higher
dimensions the relation between $\,\delta'$ and $\,\delta'_s\,$ is
not that simple.}
   \end{remark}

We shall use again the definitions (\ref{parameters}) and employ the
symbol $\,\sigma(\ell,\theta,D)\,$ for the spectrum of the
$\,\delta'_s$ lattice Hamiltonian. The conditions (\ref{delta'_s
band}) and (\ref{delta'_s band-}) have the following easy
consequences:

   \begin{proposition} \label{delta'_s basic}
$\;$ (a) The spectrum has a band structure. For any nonzero $\,D\,$
the number of gaps is infinite, $\,\sigma(\ell,\theta,D)=
\bigcup_{r=1}^{\infty} [\alpha_r,\beta_r]\,$ with $\,\alpha_r<
\beta_r<\alpha_{r+1}\,$.
\vspace{.8em}

\noindent
(b) $\;\sigma(\ell,\theta,0)= [0,\infty)\,$.
\vspace{.8em}

\noindent
(c) If $\,D>0\,$, each $\,\alpha_r\,$ equals $\,\left({\pi
n\over\lambda} \right)^2\,$ for $\,\lambda=\ell_1\,$ or $\,\ell_2\,$,
and some $\,n\in\Z\,$. On the other hand, $\,\beta_r=\left({\pi
m\over\lambda} \right)^2\,$ with $\,m\in\Z\,$ for $\,D<0\,$ and
$\,r\ge 2\,$.
\vspace{.8em}

\noindent
(d) $\;\alpha_1=0\;$ for $\;D\ge 0\,$, while $\;\alpha_1<0\;$ if
$\;D<0\,$.
\vspace{.8em}

\noindent
(e) For $\,-\ell_1\!-\ell_2<D<0\,$ we have $\,\beta_1<0\,$ and
$\;\alpha_2=0\,$.
\vspace{.8em}

\noindent
(f) $\,\sigma(\ell,\theta,D')\cap\R_+ \,\subset\,
\sigma(\ell,\theta,D)\cap\R_+\,$ holds for $\,|D'|>|D|\,$.
\vspace{.8em}

\noindent
(g) Each gap is contained in the intersection of a pair of
$\,\delta'\,$ Kronig--Penney gaps --- see Remark~\ref{delta'
rem} --- with the coupling constant $\,D\,$ and spacings $\,\ell_1\,$
and $\,\ell_2\,$.
   \end{proposition}

In distinction to the $\,\delta\,$ lattice there are therefore
always infinitely many gaps for $\,D\ne 0\,$. As in the
Kronig--Penney case and its $\,\delta'$ analogue, the roles of bands
and gaps are, roughly speaking, reversed. Comparing to the part (h)
of Proposition~\ref{delta basic}, however, the asymptotics of band
widths is slightly more complicated.

It is clear that if a band $\,[\alpha_r,\beta_r]\,$ with $\,r\gg 1\,$
is well separated from the rest of the spectrum, its width has the
same leading term as in the Kronig--Penney situation,
   \begin{equation} \label{single band}
\beta_r-\,\alpha_r \,=\, {8\over D\ell_j}\,+\,\OO(r^{-1})\,,
   \end{equation}
where $\,\ell_j\,$ is the length to which this band corresponds by
Proposition~\ref{delta basic}c. If $\,\theta\,$ is rational and the
points $\,k_n\,$ and $\,\tilde k_m\,$ coincide, we get in the same
way
   \begin{equation} \label{double band}
\beta_r-\alpha_r\,=\, {8\over D}\,(\ell_1^{-1}\!+\! \ell_2^{-1})
+\,\OO(r^{-1})\,.
   \end{equation}
It may happen, however, that $\,k_n\,$ and $\,\tilde k_m\,$ are not
identical but close to each other, so that they still produce a
single band. It is obvious from (\ref{delta'_s band}) that this leads
to an enhancement of the band width above the value given by
(\ref{double band}). The effect is most profound just before the band
splits.  Suppose, \eg, that $\,D>0\,$ and $\,\tilde k_m>k_n\,$ with
$\,\delta:= \tilde k_m\!-k_n \approx 8/D\ell_1\,$. The band width in
the momentum variable is then $\,\delta\!+\!\eta\,$ up to error
terms, where $\,\eta\,$ solves the equation
$$
\eta^2- \delta\theta^{-1}\eta- \delta^2\theta^{-1}\,=\,0\,.
$$
The enhancement due to band conspiracy is therefore
   \begin{equation} \label{enhancement}
{{\delta+\eta}\over{\delta(1+\theta^{-1})}} \,=\, g(\theta) \,:=\,
{{2\theta+1+ \sqrt{1+4\theta}}\over {2(\theta+1)}}\,.
   \end{equation}
The same can be done if the order of $\,k_n,\,\tilde k_m\,$ is
reversed; in combination with (\ref{single band}) and (\ref{double
band}) this yields the following estimates on the band width,
   \begin{equation} \label{band width bound}
{8\over DL}\,+\,\OO(r^{-1}) \,<\, \beta_r-\alpha_r\,<\, {8\over
D}\,(\ell_1^{-1}\!+\! \ell_2^{-1})\,e(\theta)\,+\,\OO(r^{-1})\,,
   \end{equation}
where
$$
e(\theta)\,:=\, \max\{ g(\theta),g(\theta^{-1})\}\,.
$$
It is easy to see that $\,e(\theta)>1\,$ with $\,\lim_{\theta\to 0}
e(\theta)= \lim_{\theta\to\infty} e(\theta)= 1\,$. For $\,\theta=1\,$
we have $\,e(\theta)=\,{1\over 4}\,(3\!+\!\sqrt 5)= 1.309\dots\,$,
but the strongest conspiracy occurs if the lattice spacing is close
to two to one, and the wider of the conspiring bands is below the
other one (above for $\,D<0\,$), because the right side of (\ref{band
width bound}) reaches its maximum at $\,g(2)=\,{4\over 3}\,$. The
value $\,\theta=2\,$ itself is, of course, integer so there is no
enhancement. Summing this discussion, we have

   \begin{proposition} \label{delta'_s band width}
The band widths of a $\,\delta'_s\,$ lattice with $\,D\ne 0\,$
satisfy the asymptotic bounds (\ref{band width bound}), where
$\,e(\theta)< \,{4\over 3}\,$.
   \end{proposition}

\setcounter{equation}{0}
\section{Conclusions}

We have said in the introduction that a choice of the coupling at
graph nodes in a realistic model should follow from a suitable
``zero--diameter'' limit of a tube system of the same topology.
To illustrate some problems which may arise, recall that a
cross--type region in the plane with the Dirichlet boundary exhibits
a bound state \cite{SRW}. Moreover, using the Dirichlet bracketing
\cite[Sec.XIII.15]{RS} in combination with the results of \cite{ABG},
one can check easily that {\em any} branched (star--shaped) system of
infinitely long tubes with the Dirichlet boundary has at least one
bound state \cite{ES2}; sometimes the number of bound states may be
even large,\eg, for systems of a narrow--fork form as can be seen
from \cite{ABG}.  This conclusion extends to situations where the
connecting--region boundary is ``rounded'' (no squeezing allowed) and
the tubes involved are only {\em asymptotically} straight
\cite{DE,ES2,ES3}.

The existence of a single bound state could be preserved in the
zero--diameter limit provided the corresponding ``coupling constant''
is chosen with a proper sign as we have pointed out in Section~2.2.
However, the couplings discussed here cannot accomodate multiple
bound states.  Using the more general boundary conditions
(\ref{permutation bc}) does not help: it is easy to see that such a
coupling has at most {\em two} bound states,
$\,\epsilon_j=-\kappa_j^2\,$ with
   \begin{equation} \label{permutation bound states}
\kappa_1\,:=\, -\,{1\over A-B}\,,\qquad \kappa_2\,:=\, -\,{1\over
A+(n-1)B}\;;
   \end{equation}
this happens if the denominators are negative. In general therefore a
zero--diameter limit should be expected to work in an energy interval
around the continuous spectrum threshold which is kept fixed, or at
most it remains small with respect to the intermode distances when
the junction region is scaled.

The analysis of both a single junction and a rectangular lattice
shows the exceptional role of the ``free'' operators, \ie, those
having the $\,\delta\,$ coupling with $\,c=0\;$ (which is the same as
$\,\delta'$ with $\,C=0\,$), or $\,\delta'_s\,$ coupling with
$\,D=0\,$. Their S--matrix elements, band profile on a lattice, and
other properties are, of course, nontrivial due the branching;
however, they are the simplest possible. For instance, the reflection
and transmission amplitudes through an $\,n\,$ wire junction are
   \begin{equation} \label{free rt}
r\,=\,\pm\,{2-n\over n}\,,\qquad t\,=\,\pm\,{2\over n}
   \end{equation}
for the ``free'' $\,\delta\,$ and $\,\delta'_s\,$ coupling,
respectively, independently of energy. Numerical calculations of
transport properties through Y--junctions \cite{Me} suggest that at
least for some systems of $\,n\,$ coupled straight tubes, this might
be the correct low--energy scattering limit.

On the other hand, a junction of finite--width tubes can have various
geometries; in fact, an experimentalist would hardly guarantee that
three joined quantum wires have a perfect Y shape. Moreover, if the
connection region supports a potential, albeit a weak one, the
low--energy scattering properties would be substantially altered.
Hence the ``non--free'' boundary conditions are also of physical
interest; one can even conceive easily a tube system in which
junction parameters are tuned by application of an external field.

Up to now we had in mind mostly simple junctions. Replacing them by
regions of a nontrivial topology we arrive at a situation to which
the considerations of Section~2.3 might be regarded as a simplified
picture.  Since we have shown there that the $\,\delta'$ and
$\,\delta'_s\,$ couplings are, at least within a fixed interval of
``intermediate'' energies, modelled by complicated enough geometric
scatterers, also the conditions (\ref{C}) and (\ref{D}) are likely to
have something in common with the real world. Moreover, we have seen
that the coupling constant of the limiting ideal scatterer is nonzero
and it is fully determined by the geometrical properties of its
approximants.

In this respect, a comment is due. Without giving any details, the
authors of \cite{GLR} suggested that such ``composed'' junctions can
be described by the boundary conditions (\ref{delta bc}) with the
``renormalized'' parameter dependent on energy. This does not
contradict our conclusions. For instance, in the simplest case
$\,n=2\,$ the corresponding reflection and transmission amplitudes,
which differ just by the sign from (\ref{rt's}), may be written
formally as the corresponding $\,\delta\,$ scattering quantities
provided we choose $\,c(k):=-Dk^2\,$. Hence if one wants to describe
a ``composed'' junction by means of a ``dressed'' coupling constant
--- which anyhow makes sense only when a prescription to compute the
latter is given --- it may happen that it differs substantially from
the ``bare'' coupling.

With this we leave this subject and turn to latice Hamiltonians of
Section~3. In addition to their possible use as models of quantum
wire superlattices, they represent an interesting mathematical
object, and the observed dependence of the spectra on
number--theoretical properties of the parameter $\,\theta\,$ raises
many questions. One would like to know, for example, how the band and
gap patterns do actually look, what are their fractal properties with
respect to $\,\theta\,$, or what is the measure of the spectrum
relative to a suitable measure on $\,\R\,$. We intend to
return to these problems in a later publication.

The results of Sections~3.2 and 3.3 show that despite there are
``less'' gaps in two dimensions, and despite the behaviour of bands
and gaps is in general irregular --- for an irrational $\,\theta\,$
they exhibit asymptotically a ``squared quasiperiodic'' distribution
--- it coincides roughly with that of the Kronig--Penney analogues to
our lattices, namely that for the $\,\delta\,$ coupling the bands
dominate at high energies, while the converse is true for the
$\,\delta'_s\,$. In analogy with the one--dimensional case
\cite{AEL1,AEL2,MS} , one can therefore make a conjecture concerning
the situation when a $\,\delta'_s\,$ lattice with $\,D\ne 0\,$ is
placed into an electric field. The heuristic tilted--band picture
suggests the existence of localization; the spectrum will remain
continuous, of course, but an unrestricted propagation may be
possible only in the direction {\em perpendicular} to the electric
field. In the $\,\delta\,$ lattice case, where we have for the
one--dimensional situation a guess but no rigorous result, the
problem is even more exciting; the results of the present paper show
that at least for some values of the lattice parameters there is no
localization.
\vspace{5mm}

\noindent
{\bf Acknowledgment.} The author is grateful for the hospitality
extended to him in the Institute of Mathematics, University of Ruhr,
Bochum, where this work was done. It is a pleasure to thank Y.~Last
for some helpful remarks. The research has been partially supported
by the Grant AS No.148409 and the European Union Project
ERB--CiPA--3510--CT--920704/704.

\vspace{30mm}
\section*{Figure captions}
\vspace{5mm}

\noindent
{\bf Figure 1.} Scattering on $\,\Gamma_n(N,\ell)\,$ with $\,n=3\,$
and $\,N=2\,$.
\vspace{5mm}

\noindent
{\bf Figure 2.} A rectangular latice.

\end{document}